%% file: CactusGraphandGCC.tex
\documentclass[11pt]{article}
\usepackage{ graphics}

\newcommand{\mathsym}[1]{{}}
\newcommand{\unicode}[1]{{}}
\usepackage{latexsym,xy,eucal,mathrsfs,graphs}
\usepackage{setspace}
\usepackage{amsthm}
\usepackage{pstricks}
\usepackage{amssymb}
\usepackage{amsmath}
\usepackage{indentfirst}
\usepackage{epsfig}
\usepackage{indentfirst}
\usepackage{amscd}
\usepackage{multirow}

\newtheorem{theorem}{Theorem}[section]
\newtheorem{proposition}[theorem]{Proposition}
\newtheorem{lemma}[theorem]{Lemma}
\newtheorem{claim}[theorem]{CLAIM}

\newtheorem{corollary}[theorem]{Corollary}

\newtheorem{result}[theorem]{Result}
\theoremstyle{definition}
\newtheorem{defn}[theorem]{Definition}
\newtheorem{example}[theorem]{Example}

\usepackage{simplemargins}
\setleftmargin{0.88in}
\setrightmargin{0.88in}
\settopmargin{0.817in}
\setbottommargin{0.83in}

 \font\conj =
msbm10 at 12pt  

\def\text#1{\mathop{{\rm{#1}}}\nolimits}
\def\mr{\mathop{{\rm{mr}}}\nolimits}
\def\msr{\mathop{{\rm{msr}}}\nolimits}
\def\max{\mathop{{\rm{max}}}\nolimits}

\def\nulli{\mathop{{\rm{null}}}\nolimits}

\def\cc{\mathop{{\rm{cc}}}\nolimits}

\def\d{\mathop{{\rm{d}}}\nolimits}

\def\gcc{\mathop{{\rm{GCC}}}\nolimits}
\def\gccp{\mathop{{\rm{GCC}_+}}\nolimits}

\def\Re{{\mbox{\conj R}}}

\def\F{{\mbox{\conj F}}}
\def\C{{\mbox{\conj C}}}

\def\dis{\displaystyle}
\def\1{\'{\i}}

\title{Cactus Graphs  and Graph Complement  Conjecture}

\author{Pedro D\'{\i}az Navarro\thanks{Escuela de Matem\'atica, Universidad de Costa Rica}}
\date{August, 2016}
\begin{document}
\maketitle

\begin{abstract}
\noindent
In  this paper   we  proof  that  any  cactus  graph   satisfies   graph  complement  conjecture  by  finding  a orthogonal  representation  of  its  complement  in  $\Re^5$.
\end{abstract}

\\

\noindent{\bf DOI:} 05C50,05C76 ,05C85 ,68R05 ,65F99,97K30

\section{Introduction}

A {\it graph} $G$ consists of a set  of vertices $V(G)=\{1,2,\dots,n\}$ and a set of edges $E(G)$, where an edge is defined to be an unordered  pair of vertices. The {\it order} of $G$, denoted $\vert G\vert $,  is the cardinality  of $V(G)$. A graph is {\it simple} if it has no multiple  edges  or loops. The {\it complement } of a graph $G(V,E)$ is the graph $\overline{G}=(V,\overline{E})$, where $\overline{E}$ consists of all those edges of the complete  graph $K_{\vert G\vert}$ that are not in $E$.

 A matrix $A=[a_{ij}]$ is {\it combinatorially symmetric} when $a_{ij}=0$ if and only if $a_{ji}=0$. We say that  $G(A)$  is the graph of a combinatorially symmetric matrix $A=[a_{ij}]$ if $V=\{1,2,\dots,n\}$ and $E=\{\{i,j\}: a_{ij}\ne0\}$ . The main diagonal entries of $A$ play no role in determining $G$. Define $S(G,\F)$ as the set of all $n\times n$ matrices that are real symmetric if $\F=\Re$ or complex Hermitian if $\F=\C$ whose graph is $G$. The sets $S_+(G,\F)$ are the corresponding subsets of positive semidefinite (psd) matrices. The smallest possible rank  of any matrix $A\in S(G,\F)$  is the {\it minimum rank} of $G$, denoted $\mr(G,\F)$, and the smallest possible rank of any matrix $A\in S_+(G,\F)$  is the {\it minimum semidefinite rank} of $G$, denoted $\mr_+(G)$ or $\msr(G)$.

In 1996, the minimum rank among real symmetric matrices with a given graph was studied  by  Nylen \cite{PN}. It gave rise to the area of minimum rank problems which led to the study of minimum rank among complex Hermitian matrices and positive semidefinite matrices associated with a given graph. Many results can be  found for  example  in  \cite{FW2, VH, YL, LM, PN}.

During  the   AIM workshop of 2006 in Palo Alto, CA, it  was  proposed  question about  how large can  $\mr(G) + \mr(\overline{G})$ be \cite{RB1}the     that for any   graph $G$  and  infinite field  $F$.  It was conjectured that $\mr(G)+ \mr(\overline{G})\le \vert G\vert +2$   for some   infinite  families of graphs for which the minimum rank of both  the graph and their complement were known  but   for an arbitrary  graph $G$  it is unknown whether or not   this inequality is true and therefore it  is called the ``{\it Graph Complement Conjecture}'' denoted  as $\gcc$. When  we  restrict  the study  to $S_+(G,\F),  F=\C\ \hbox{or}\ \Re$  this  conjecture is refered to as $\gccp$,  which  states  that  for  any graph $G$, $\msr(G)+\msr(\overline{G})\le |G|+2$.

In \cite{FW2} it is shown  that  trees satisfy  $\hbox{\rm GCC}_+$ since the minimum semidefinite  rank of a tree is $\vert G\vert -1$ and that of its complement was shown to be  at most 3. Hogben \cite{LH2},  showed that   some other families of sparse graphs, including unicyclic graphs and $2$-trees also satisfy  $\hbox{\rm GCC}_+$. Sharawi \cite{SY}, showed that $\hbox{\rm GCC}_+$  holds for a complete bipartite graph $G$ with $V(G)= R\sqcup L, \vert L\vert=n\ge 2, \vert R\vert = m\ge n$.  Also,  it was shown in \cite{SY} that  a $k-$regular Harary graph, $H_{k,n}$ which is $k$-connected satisfies  the $\hbox{\rm GCC}_+$. Mitchell \cite{LM},  showed  that chordal graphs satisfy $\hbox{\rm GCC}_+$.  Other results   have  been  found in  \cite{FW, JE, LH2, LM, SY}. However,  the  general problem  of    graph  complement  conjecture is still open.

\section{Graph Theory Preliminaries}
\addtocontents{toc}{\vspace{-15pt}}
In  this section   we give   definitions and results from  graph theory which    will  be used in  the remaining  sections. Further details  can be found  in \cite{BO,BM, CH}.

A {\bf graph} {$G(V,E)$} is a pair {$(V(G),E(G)),$} where {$V(G)$} is the set of vertices and {$E(G)$} is the set of edges together  with an  {\bf  incidence  function} $\psi(G)$ that associate with  each edge  of  $G$ an  unordered  pair  of (not necessarily  distinct) vertices  of  $G$. The {\bf order} of {$G$}, denoted {$|G|$}, is the number of vertices in {$G.$} A graph is said to be {\bf simple} if it has no loops or multiple edges. The {\bf complement} of a graph {$G(V,E)$} is the graph {$\overline{G}=(V,\overline{E}),$} where {$\overline{E}$} consists of all the edges that are not in {$E$}.
 A {\bf  subgraph} {$H=(V(H),E(H))$} of {$G=(V,E)$} is a graph with {$V(H)\subseteq V(G)$} and {$E(H)\subseteq E(G).$} An {\bf induced subgraph} {$H$} of {$G$}, denoted G[V(H)], is a subgraph with {$V(H)\subseteq V(G)$} and {$E(H)=\{\{i,j\} \in E(G):i,j\in V(H)\}$}. Sometimes  we  denote the  edge $\{i,j\}$ as $ij$.
 We  say  that  two  vertices of a graph $G$  are {\bf adjacent}, denoted  $v_i\sim v_j$,   if  there is an edge $\{v_i,v_j\}$ in  $G$.  Otherwise  we say  that the  two  vertices $v_i$ and $v_j$ are {\bf non-adjacent}  and  we denote this  by $v_i \not\sim v_j$.  Let {$N(v)$} denote the set of vertices that are adjacent to the vertex {$v$} and let {$N[v]=\{v\}\cup N(v)$}. The {\bf degree} of a vertex {$v$} in {$G,$} denoted {$\d_G(v),$} is the cardinality of {$N(v).$} If {$\d_G(v)=1,$} then {$v$} is said to be a {\bf pendant} vertex of {$G.$} We use {$\delta(G)$} to denote the minimum degree of the vertices in {$G$}, whereas {$\Delta(G)$} will denote the maximum degree of the vertices in {$G$}.
 Two  graphs $G(V,E)$ and $H(V',E')$  are  identical  denoted  $G=H$, if  $V=V',  E=E'$, and $\psi_G=\psi_H$  . Two  graphs $G(V,E)$ and $H(V',E')$ are {\bf isomorphic}, denoted  by  $G\cong H$, if  there exist bijections $\theta:V\to V'$  and $\phi: E\to  E' $ such  that $\psi_G(e)=\{u,v\}$  if  and  only if  $\psi_H(\phi(e))= \{\theta(u), \theta(v)\}$.
 A {\bf complete graph}  is a simple graph in which the vertices are pairwise adjacent.
  We will use {$nG$} to denote {$n$} copies of a graph {$G$}. For example, $3K_1$  denotes three  isolated vertices $K_1$ while {$2K_2$} is the graph given  by   two  disconnected  copies  of $K_2$.
 A {\bf path} is a list of distinct vertices in which successive vertices are connected by edges. A path on {$n$} vertices is denoted by {$P_n.$} A graph {$G$} is said to be {\bf connected} if there is a path between any two vertices of {$G$}. A {\bf cycle} on {$n$} vertices, denoted {$C_n,$} is a path such that the beginning vertex and the end vertex are the same. A {\bf tree} is a connected graph with no cycles. A graph $G(V,E)$ is  said  to be {\bf chordal}  if it has no induced cycles $C_n$ with $n\ge 4$.
 A  {\bf component}  of a graph $G(V,E)$ is  a maximal connected  subgraph. A  {\bf cut vertex}   is  a vertex  whose deletion  increases  the number of  components.
 The {\bf union}  $G\cup G_2$ of  two  graphs $G_1(V_1,E_1)$ and $G_2(V_2,G_2)$  is  the union  of  their  vertex  set  and  edge  set,  that is $G\cup G_2(V_1\cup V_2,E_1\cup E_2$. When  $V_1$ and $V_2$ are disjoint their union  is called  {\bf  disjoint union} and  denoted $G_1\sqcup G_2$.


\section{The Minimum  Semidefinite Rank  of  a Graph}
In  this section  we will establish   some of  the   results for  the minimum  semidefinite rank ($\msr$)of a graph $G$  that  we  will be using in the subsequent sections.

A {\bf positive  definite} matrix  $A$ is an Hermitian   $n\times n$ matrix such  that $x^\star A x>0$  for all nonzero  $x\in \C^n$. Equivalently,  $A$  is a  $n\times n$ Hermitian positive definite matrix  if and  only  if   all the  eigenvalues of $A$ are positive (\cite{RC}, p.250).

A $n\times n$ Hermitian matrix  $A$ such  that $x^\star A x\ge 0$  for all $x\in \C^n$ is  said  to be  {\bf positive  semidefinite (psd)}. Equi\-va\-lently,   $A$ is a  $n\times n$ Hemitian positive  semidefinite matrix if and  only  if  $A$  has  all   eigenvalues nonnegative (\cite{RC}, p.182).

If $\overrightarrow{V}=\{\overrightarrow{v_1},\overrightarrow{v_2},\dots, \overrightarrow{v_n}\}\subset \Re^m$  is a set of  column vectors  then  the  matrix
$ A^T A$, where $A= \left[\begin{array}{cccc}
  \overrightarrow{v_1} & \overrightarrow{v_2} &\dots& \overrightarrow{v_n}
\end{array}\right]$
and $A^T$  represents  the  transpose matrix of  $A$, is a psd matrix  called  the  {\bf Gram matrix} of $\overrightarrow{V}$. Let $G(V,E)$  be a graph  associated  with  this Gram matrix. Then  $V_G=\{v_1,\dots, v_n\}$ correspond to  the set of  vectors in $\overrightarrow{V}$ and  E(G) correspond to  the nonzero inner products  among  the  vectors  in $\overrightarrow{V}$. In this  case $\overrightarrow{V}$  is  called an  {\bf orthogonal representation} of $G(V,E)$ in $\Re^m$. If  such  an  orthogonal  representation  exists  for  $G$ then $\msr(G)\le m$.

The  {\bf  maximum positive  semidefinite   nullity }  of a graph $G$, denoted $M_+(G)$   is   defined  by
$
M_+(G) =\max\{\nulli(A): A\ \hbox{ is symmetric and positive semidefinite and }\   G(A) =G \}
$,
where $G(A)$  is  the graph obtained  from  the matrix $A$. From  the  rank-nullity theorem we get  $ \msr(G)+ M_+(G)=|G|$.

Some  of  the most  common   results about  the minimum semidefinite  rank of a graph  are  the  following:

\begin{result}\cite{VH}\label{msrtree}
If  $T $ is a tree  then $\msr(T)= |T|-1$.
\end{result}
\begin{result}\cite{MP3}\label{msrcycle}
 The cycle $C_n$ has minimum semidefinite rank $n-2$.
\end{result}


\begin{result}\label{res2}
 \cite{MP3} \ If a connected  graph $G$  has a pendant  vertex $v$, then $\msr(G)=\msr(G-v)+1$ where $G-v$  is obtained as an induced subgraph of $G$ by  deleting $v$.
\end{result}

\begin{result} \cite{PB} \label{OS2}
 If {$G$} is a connected, chordal graph, then {$\msr(G)=\cc(G).$}
\end{result}

\begin{result}\label{res1}
\cite{MP2}\ If a graph $G(V,E)$ has a cut vertex, so that $G=G_1\cdot G_2$, then  $\msr(G)= \msr(G_1)+\msr(G_2)$.
\end{result}

The  next  two  definitions give  us   two  families of   graphs   which  are important  in  the study   of  the minimum semidefinite  ranbk of  a  graph.
\begin{defn}\label{ccpg}
Suppose  that $G=(V,E)$  with $|G|=n \ge 4$  is simple  and  connected such  that  $\overline{G}=(V,\overline{E})$  is also  simple  and  connected. We  say  that  $G$ is a {\bf $\mathbf{\delta}$-graph} if  we  can  label  the vertices of $G$ in such a way that
\begin{enumerate}
  \item[(1)] the  induced    graph    of  the  vertices  $v_1,v_2,v_3$  in $G$ is  either $3K_1$  or  $K_2 \sqcup K_1$,  and
  \item[(2)] for $m\ge 4$,  the vertex  $v_m$  is  adjacent   to all   the  prior  vertices $v_1,v_2,\dots,v_{m-1}$  except  for at most $\dis{\left\lfloor\frac{m}{2}-1\right\rfloor}$ vertices.
  \end{enumerate}
 \end{defn}
\begin{defn} Suppose   that  a graph $G(V,E)$  with $|G|=n \ge 4$  is  simple and  connected   such  that  $\overline{G}=(V,\overline{E})$  is also   simple and  connected.  We  say  that  $G(V,E)$  is a {\bf C-$\mathbf{\delta}$  graph}  if    $\overline{G}$   is a  $\delta$-graph.

In other  words,  $G$   is a  {\bf C-$\mathbf{\delta}$ graph} if  we can  label   the  vertices   of  $G$  in  such a  way   that
 \begin{enumerate}
 \item[ (1)] the  induced  graph  of  the vertices  $v_1,v_2,v_3$ in  $G$ is  either  $K_3$ or $P_3$,  and
 \item[(2)]  for  $m\ge 4$,  the vertex $v_m$  is adjacent  to at most  $\dis{\left\lfloor\frac{m}{2}-1\right\rfloor}$ of  the  prior  vertices $v_1,v_2,\dots,v_{m-1}$.
\end{enumerate}
\end{defn}
\begin{example}\label{examplecp}
The   cycle $C_n, n\ge 6$  is a C-$\delta$ graph and  its  complement  is  a  $\delta$-graph.
\begin{center}
\includegraphics[height=30mm]{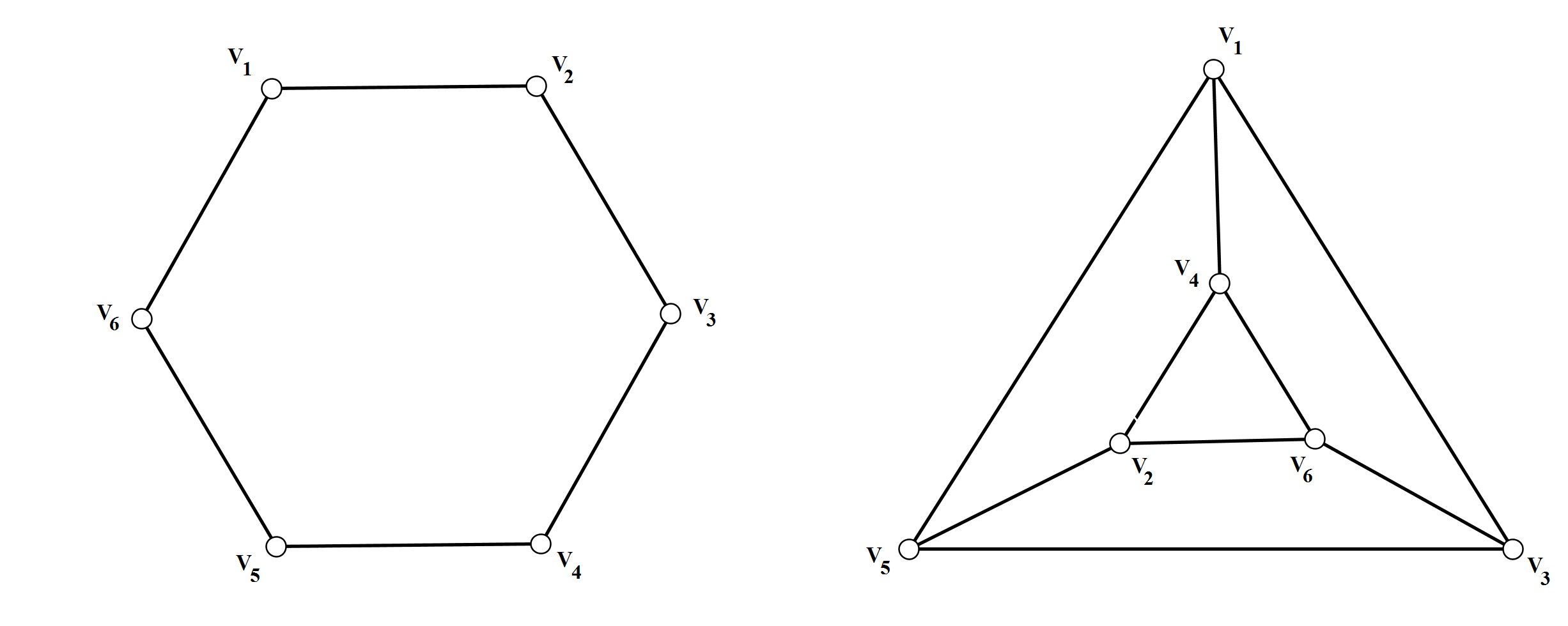}
 \end{center}
\vspace{-0.1in}\begin{figure}[h]
\centering
\caption{The  Graphs $C_6$  and  the  3-Prism }\label{fig3.1}
\end{figure}
Note that we  can label  the vertices  of $C_6$ clockwise $v_1,v_2,v_3,v_4,v_5,v_6$. The  graph  induced   by  $v_1,v_2,v_3$  is $P_3$. The    vertex $v_4$ is   adjacent to a prior vertex   which  is  $v_3$. Also,   the  vertex $v_5$ is adjacent to  vertex $v_4$ and   the  vertex $v_6$ is adjacent to two prior  vertices  $v_1$ and $v_5$.  Hence, $C_6$  is C-$\delta$ graph. The $3$-prism which  is  isomorphic to  the complement  of  $C_6$,  is a $\delta$-graph.
\end{example}
 In  \cite{PD}it was proved  the  following  result which  prove  that  any  $\delta$-graph satisfies  delta  conjecture.
\begin{theorem} \label{main} Let $G(V,E)$ be a $\delta$-graph then
$$
\msr(G)\le\Delta(\overline{G})+1=|G|-\delta(G)\label{mrsineq1}
$$
\end{theorem}

\section{Cactus Graphs and Graph  Complement  Conjecture}

In  this  section  we  will  prove  that the complement  of   any  cactus  graph  has an  orthogonal representation  in $\Re^5$.

\begin{defn}\cite{BLS}
 A simple connected graph $G$ is a {\bf cactus  graph} if every edge is part of at most one cycle in $G$.
\end{defn}
\begin{center}
\includegraphics[height=45mm]{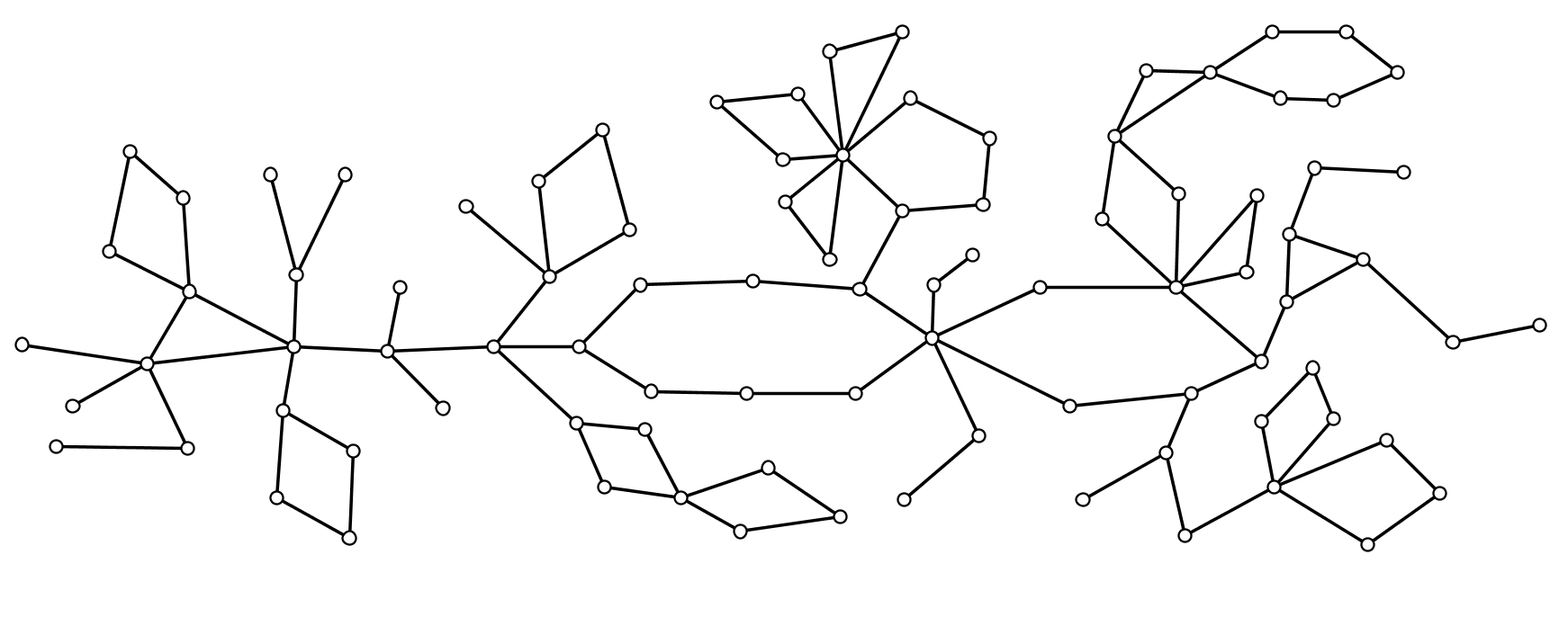}
\end{center}
\vspace{-0.1in}\begin{figure}[h]
\centering
\caption{Cactus Graph } \label{cactus1}
\end{figure}

Equivalently, a  simple connected graph $G$  is a cactus  graph  if   any two simple cycles in $G$ have at most one vertex in common. From \cite{BLS} we know that a  cactus  graph $G$ is outerplanar  since  $G$  cannot  contain  $K_4$ or $K_{2,3}$  as  a minor.
\begin{example} Some familiar  simple connected  graphs  are  cactus  graphs. For  example:
\begin{enumerate}
\item Trees    are cactus graphs  with no cycles.
\item  Unicyclic  graphs   are  cactus  graphs  with  only  one  cycle.
\item Chains  are   cactus   graphs   with at least  two cycles.
\end{enumerate}
\end{example}
 It  is possible  to label  the vertices  of   a  cactus  graph $G, |G|> 5$, satisfying the definition  of a C-$\delta$ graph.  As a consequence, the  graph complement of   any  cactus  graph   has an  orthogonal   representation  in  $\Re^{\Delta (G)+1}$. However, since  the number  of  cycles  sharing a single  vertex  in  a cactus  graph  could  be arbitrarily large, the upper  bound $\Delta (G)+1$  for  $\msr(\overline{G})$ is   too large   to   prove $\gccp$. The  next  proposition   gives   sufficient  conditions    to get  an  orthogonal  representation   in $\Re^5$  for  the  graph  complement  of  a simple connected cactus graph.    We  will use  this  result to  prove    that  any  cactus  graph  satisfies $\gccp$.
\begin{proposition}\label{proptcg}
Let $G(V,E), |G|\ge 5$  be a simple connected  graph  that  can be  constructed  from  a  path  $P: v_1v_2v_3$  in  such a way   that  the  newly added vertex $v_m, m\ge 5$  is adjacent to at most  two of  the prior  vertices $v_1,v_2,\dots,v_{m-1}$. Then  $\overline{G}(V,\overline{E})$  has an orthogonal representation of pairwise linearly independent vectors in $\Re^5$.
\end{proposition}
\noindent\MakeUppercase{Proof}:

If $G$  is a simple connected   graph   with  $|G|\le 5$  it is  easy  to  check  that $\overline{G}$  has an  orthogonal representation  in $\Re^5$.  So  Assume  that $|G|=n>5$. Let $V= \{v_1,v_2,\dots,v_n\}$. There is a path $P_3: v_1v_2v_3$  as a subgraph of  $G$ induced   by $\{v_1,v_2,v_3\}\subseteq V_G$. Let $H$ be the induced  graph of $\overline{G}$ obtained   from  the  vertices $v_1,v_2$ and $v_3$.  Then $H$ is  $K_2\sqcup K_1$. Let $\{\overrightarrow{e}_1,\overrightarrow{e}_2,\overrightarrow{e}_3,\overrightarrow{e}_4,\overrightarrow{e}_5\}$ be the  standard  orthonormal  basis  for $\Re^5$.

Using  a  similar argument  as  in Theorem  \ref{main}   we can  get  vectors  $\overrightarrow{v}_1,\overrightarrow{v}_2,\overrightarrow{v}_3$ and $\overrightarrow{v}_4$ in  $\Re^5$ for  the  vertices $v_1,v_2,v_3$, and $v_4$  where  all  the  entries  of  these vectors  are nonzero  and  belong  to  different  field  extensions.

Assume  that  for  any,  $Y_{m-1}=(V_{Y_{m-1}},\overline{E}_{Y_{m-1}}), V_{Y_{m-1}}= \{v_1,v_2,\dots,v_{m-1}\},5\le m\le n$
it is possible  to get  an orthogonal representation of pairwise linearly independent vectors $\overrightarrow{v}_1,\overrightarrow{v}_2,\dots,\\ \overrightarrow{v}_{m-1}$ in $\Re^5$, of  the form
$
 \overrightarrow{v}_i= k_{i,1}\overrightarrow{e}_1+k_{i,2}\overrightarrow{e}_2+k_{i,3}\overrightarrow{e}_3+k_{i,4}\overrightarrow{e}_4+k_{i,5}\overrightarrow{e}_5,  i= 1,\dots,m-1
$\
such  that  all  the entries of  vectors  can be  chosen  nonzero and   from  different  field extensions. Let $ Y_m$   be  the induced  graph  of  $G$   given  by   the  vertices  $v_1,v_2,\dots, v_m$.
 Assume  that  $v_m$  has  a  vector
$
\overrightarrow{v}_m= k_{m,1}\overrightarrow{e_1}+ k_{m,2}\overrightarrow{e_2}+\dots+ k_{m, 5}\overrightarrow{e}_5.
$

\noindent{Case 1.} $v_{m}$ is adjacent in  $G$ to two  vertices $v_i$ and $v_j, i\ne j$  from $\{v_1,v_2,\dots, v_{m-1}\}$.

 Let $\rho$ be a permutation of $(1,2,\dots,m-1)$. Suppose $v_{\rho(1)},v_{\rho(2)},\dots,v_{\rho(m-3)}$ are  adjacent to $v_m$ in $\overline{G}$  and $v_{\rho(m-2)},
v_{\rho(m-1)}$  are  not  adjacent   to  $v_{m}$ in  $\overline{G}$. The vectors  $\overrightarrow{v}_{\rho(1)},\overrightarrow{v}_{\rho(2)},\dots,\overrightarrow{v}_{\rho(m-3)},\\  \overrightarrow{v}_{\rho(m-2)}, \overrightarrow{v}_{\rho(m-1)}$ and $\overrightarrow{v}_{m}$  satisfy  the non-homogeneous system  $S$ given by:
\begin{eqnarray*}
  \langle\overrightarrow{v}_{\rho(1)},\overrightarrow{v}_m\rangle &=&g_{m,1},\ \  g_{m,1}\ne 0\\
 \langle\overrightarrow{v}_{\rho(2)},\overrightarrow{v}_m\rangle &=&g_{m,2},\ \  g_{m,2}\ne 0\\
 \vdots&\vdots&\vdots\\
 \langle\overrightarrow{v}_{\rho(m-3)},\overrightarrow{v}_m\rangle &=&g_{m,m-3},\ \  g_{m,m-3}\ne 0\\
\langle\overrightarrow{v}_{\rho(m-2)},\overrightarrow{v}_m\rangle&=& 0 \\
\langle\overrightarrow{v}_{\rho(m-1)},\overrightarrow{v}_m\rangle &=&0
\end{eqnarray*}
\noindent containing $m-3$  equations from  the  adjacency conditions in  $\overline{G}$ and  two  equations from  the  orthogonal conditions in $\overline{G}$. The vectors $\overrightarrow{v}_{\rho(i)}, i=1,2,\dots,m-1$ have  the form
  $ \overrightarrow{v}_{\rho(i)}= k_{\rho(i),1}\overrightarrow{e}_1+k_{\rho(i),2}\overrightarrow{e}_2+\dots+k_{\rho(i),5}\overrightarrow{e}_5$
where all  $k_{\rho(i),j}, i=1,2,\dots,m-1 ,  j= 1,2,\dots,5$   are not  zero.

  Similar argument  as  in the  proof  for  $\delta $-graphs  of  the  theorem   \ref{main}  can be applied  to   get  a  non-homogeneous system  $S$  in  the variables $k_{m,1},\dots, k_{m,5}$  and  a  homogeneous  system $S_H$ in the variables $k_{m,1},\dots, k_{m,5},-g_{m,1}, \dots, -g_{m,m-3}$. The matrix $A$ of  the  homogeneous  system  $S_H$  is  given  by

  {\tiny
  $$A=
\left(\begin{array}{cccccccccccc}
   k_{\rho(1),1}& k_{\rho(1),2} & k_{\rho(1),3} &k_{\rho(1),4} & k_{\rho(1),5} & -1 & 0 & 0& 0 &0 &\dots&0\\
   k_{\rho(2),1} &k_{\rho(2),2} &k_{\rho(2),3} & k_{\rho(2),4} & k_{\rho(2),5} &0 &-1 & 0&0& 0&\dots&0 \\
   k_{\rho(3),1} & k_{\rho(3),2} & k_{\rho(3),3} & k_{\rho(3),4}&k_{\rho(3),5} &0 & 0 & -1 &0 &0 &\dots&0\\
   k_{\rho(4),1} & k_{\rho(4),2} & k_{\rho(4),3} & k_{\rho(4),4}&k_{\rho(4),5} &0 & 0 & 0 &-1 & 0&\dots&0\\
   k_{\rho(4),1} & k_{\rho(4),2} & k_{\rho(4),3} & k_{\rho(4),4}&k_{\rho(4),5} &0 & 0 & 0 &0 & -1&\dots&0\\
\vdots& \vdots & \vdots & \vdots &\vdots & \vdots & \vdots &\vdots &\vdots& \vdots&\dots&\vdots \\
 k_{\rho(m-3),1} & k_{\rho(m-3),2} & k_{\rho(m-3),3} & k_{\rho(m-3),4} & k_{\rho(m-3),5} & 0 & 0 &0 & 0 &0&\dots&-1 \\
 k_{\rho(m-2),1} & k_{\rho(m-2),2}& k_{\rho(m-2),3} & k_{\rho(m-2),4}& k_{\rho(m-2),5} &0 & 0 &0 &0 & 0&\dots&0 \\
k_{\rho(m-1),1} & k_{\rho(m-1),2}& k_{\rho(m-1),3} &  k_{\rho(m-1),4}& k_{\rho(m-1),5} &0 & 0 &0 &0 & 0&\dots&0 \\
\end{array}\right)
$$
}
Reducing  this matrix  to echelon  form   we  get  the  matrix  $B$
{\tiny
$$B=
\left(\begin{array}{cccccccccccccccccccc}
  1&\ast & \ast &\ast&\ast&| & \delta_1 & 0 & 0& 0 &0&0&0&\dots& 0& 0 &0&0&0&0\\
0 &1 &\ast & \ast&\ast&| &\ast & \delta_2&0&0& 0&0&0&\dots& 0& 0 &0&0&0&0 \\
0 & 0 & 1&\ast &\ast&| &\ast& \ast& \delta_3&0&0&0&0&\dots& 0& 0 &0&0&0&0\\
0 & 0 &0 &1&\ast&| &\ast &\ast&\ast &\delta_{4}&0 &0&0&\dots& 0& 0 &0&0&0&0 \\
0 &0& 0 &0&1&| &\ast &\ast &\ast& \ast& \delta_{5}&0&0&\dots& 0& 0 &0&0&0&0  \\
-&-&-&-&-&-&-&-&-&-&-&-&-&\dots&-&-&-&-&-&- \\
0 & 0&0 &0& 0&| &1& \ast &\ast &\ast & \ast&\delta_6& 0& \dots&0&0&0&0&0&0 \\
0 & 0&0 &0& 0&| &0& 1 &\ast &\ast & \ast&\ast&\delta 7&\dots& 0 &0&0&0&0&0 \\
\vdots &\vdots&\vdots &\vdots& \vdots&| &\vdots& \vdots &\vdots& \vdots& \vdots&\vdots &\vdots& \vdots &\vdots&\vdots&\vdots&\vdots&\vdots& \vdots \\
0&0&0 &0&0&| &0& 0 &0& 0& 0&0 &0& \dots &\ast&\ast&\ast&\delta_{m-5}&0& 0 \\
0 &0&0 &0&0&| &0& 0&0& 0&0&0 &0& \dots &\ast&\ast&\ast&\ast&\delta_{m-4}& 0 \\
0 & 0 & 0& 0&0&| &0 & 0 &0&0 & 0 &0&0&\dots& 1& \ast &\ast&\ast&\ast&\delta_{m-3} \\
0 & 0 &0 &0&0&| &0 &0& 0& 0&0&0&0&\dots& 0& 1&\ast&\ast&\ast&\ast\\
0 & 0 &0 &0&0&| &0 &0& 0& 0& 0&0&0&\dots& 0& 0 &1&\ast&\ast&\ast
\end{array}\right)
$$
}
which is a block matrix
$$
B=\left(
    \begin{array}{ccc}
      B_1&| & B_2 \\
      -&&-\\
      0 &|& R \\
    \end{array}
  \right)
$$
where  $B_1$ is a square matrix  of  size $5$, the matrix $B_2$ is a matrix of  size $5\times (m-3)$,  the  zero matrix   has   size $(m-6)\times 5$, and  the  matrix  $R$ has  size $(m-6)\times (m-3)$. As a consequence, the system $R\overrightarrow{g}=0$ where  $\overrightarrow{g}=(g_{m,1},\dots, g_{m,m-3})^T$,  has  infinitely many solutions depending  on at least three free variables. Choosing  all  the free  variables   from different  field extensions as in  the  proof  of  Theorem \ref{main} we  can  get  all  the  values  $g_{m,1},\dots, g_{m,m-3}$ nonzero. Also, with a similar argument as    with  the matrix  in  Case 2 of the proof  of  \ref{main},  we  get  that   the  values   in  the diagonal  of the matrix  $B_2$ are nonzero.  As  a consequence, we   get   that $k_{m,1},\dots, k_{m,5}$ can  be chosen nonzero.

Since  the  values $g_{m,i}, i=3,4,\dots,m-3$ can be chosen  non zero and   from field extensions different  from  all previous  field extensions  taken  for  the values of $k_{i,j}, i=1,2,\dots,m, j=1,2,\dots,5$  there exist   at least  one solution  for  the system $S$. Therefore the  vector $v_{m}$  can be constructed satisfying  all  the adjacency conditions and  orthogonal conditions.

As a consequence, $Y_{m}$ has an  orthogonal  representation  of  vectors  in  $\Re^5$.  Since  $5\le m\le n=|G|$, we get that  $Y_{|G|}=\overline{G}$ has an orthogonal representation of vectors  in $\Re^5$.

\noindent{Case 2.} $v_{m}$ is adjacent in  $G$  with  only one vertex  $v_i\in \{v_1,v_2,\dots, v_{m-1}\}$.
The argument  is similar  to  the Case 1 above. Since  we  have  only one equation    from  the  orthogonal  condition $\langle \overrightarrow{v}_m,\overrightarrow{v}_i\rangle=0$ for  $\overline{G}$ we can  choose  the  values  $g_{m,1},g_{m,2},\dots,g_{m,m-3}$ and    the  values $k_{m,i}, i = 1,2,\dots,5$ nonzero.
\hfill $\Box$

Note  that  if $G$  is a graph  satisfying  the hypothesis of  Proposition \ref{proptcg}  then $\msr(\overline{G})\le 5$.

\begin{lemma}\label{lemcactus1}
Let $G(V,E)$  be a cactus  graph with  $|G|\ge 5$. Then  $\overline{G}$ can  be constructed   by  starting  with an induced   subgraph $P_3$,  by  adding  one vertex   at a time   such  that  the  newest  vertex  is adjacent   to all but at most two of  the  prior vertices.
\end{lemma}
\noindent\MakeUppercase{Proof}:

Let $G(V,E)$  be  a  cactus  graph,  $|G|\ge4$. Since   $G$ is simple  and  connected,   $G$  contains  a path  $P_3: v_1v_2v_3$.  From  Proposition \ref{proptcg}  it is  enough  to  show  that $G$  can  be constructed   from $P_3$ by  adding  one vertex   at a time   such  that  the  newest  vertex  is adjacent   to all  previous    vertices   but at most  two vertices.
Assume   that  $|G|=M$.  Starting    with  the  vertices  $v_1,v_2,v_3$  in $P_3$  and  since  $G$  is outerplanar  we  can  get  a path  $P$  traveling along   the  graph following a clockwise  orientation,  that is,  at any  intersection   we  will choose the clockwise  direction  to  continue. Then   skipping the repeated   vertices   we  can  label   all  the remaining  vertices  $v_4, v_5, \dots, v_{M-1},v_M$  following  the  order   in  which   they first appeared  in the path  $P$. Note  that   in this  way  all possible  intersection points are counted  just  once  because  after the  first  time  we arrived  at the  point  we  should  skip it. Consider  the sequence of induced  subgraphs
$Y_3(V_3,E_3)\subset Y_4(V_4,E_4)\subset\dots\subset Y_k(V_k,E_k)\subset\dots\subset Y_{M-1}(V_{M-1},E_{M-1})\subset Y_{M}(V_M,E_M) = G$
where  $V_3=\{v_1,v_2,v_3\}$  and  $V_{k} = V_{k-1} \cup \{v_{k}\}, v_{k}\not \in  Y_k$.
\begin{claim}\label{claimcactus}
  The induced   graph  of  $Y_{k}, k\ge 4$   in  $\overline{G}$ can be  constructed    in  such a way  that  $v_k$   is adjacent   to all prior  vertices $v_1,v_2,\dots,v_{k-1}$   except    to  one  or  two  vertices.
\end{claim}
\noindent\MakeUppercase{Proof of  the  Claim}:

Constructing  the  induced  graph   $Y_4$ in  $G$   we   realize  that  $v_4$ is adjacent   to $v_3$ in $G$  but  cannot  be adjacent  to more  than two   of  the  vertices $ v_1,v_2,v_3$ in  $\overline{G}$.
Assume  that  the induced  graph    $Y_{k-1}$  in  $G$   was  constructed  as  defined above. Then  we can construct  the induced  graph   $Y_k$  by  adding   the vertex $v_k$  to  the induced  graph   $Y_{k-1}$.  The  vertex  $v_k$  cannot  be one  of  the previous  vertices   because   we skip  all  the  repeated  vertices. Now   if  $v_k$  is  adjacent  to  more  than  two  of  the  prior vertices   then     the graph  $G$  contains  two  cycles   sharing  a common   edge  which  contradicts  the  definition  of  cactus  graphs.
As a consequence,  $v_k$  is adjacent  to all   but   at  most  one  or   two of  the prior  vertices  $v_1,v_2,\dots,v_{k-1}$  of  the induced  graph    $Y_{k-1}$ in  $\overline{G}$  and  the  claim  holds. Finally,  from  the  claim \ref{claimcactus},  since  $Y_M=G$  we  get  that the cactus  graph $G$  can be constructed  as stated.
\hfill $\Box$

Note  that Lemma \ref{lemcactus1} implies  that a cactus  graph $G, |G|\ge 5$  is a C-$\delta$ graph.
\begin{corollary}\label{coro1cactus}
Let   $G(V,E)$  be a cactus  graph  then $\msr(\overline{G})\le 5$.
\end{corollary}
\noindent\MakeUppercase{Proof}:

Let  $G$  be a cactus  graph. If $|G|\le 5$ it is straightforward  to  check     that $\msr(\overline{G})\le 5$.  Then assume  that $|G|\ge  5$. From  the Lemma \ref{lemcactus1} we  have  that there is  a  orthogonal representation for  $\overline{G}$  of pairwise linearly  independent  vectors in  $\Re^5$.  As a consequence  $\msr(\overline{G})\le 5$.
\hfill $\Box$
\begin{proposition}\label{prop1cactus}
Let $G(V,E)$  be a cactus  graph with  exactly  one  cycle.  Then  the   tree  cover  number $T(G)=2$.
\end{proposition}
\noindent\MakeUppercase{Proof}:

Let $G$  be a cactus  graph with  exactly  one  cycle. Then   $G$  consists of  the cycle    $C$  and   some  trees    joined  to  the cycle at some of  its  vertices.  Consider   a maximum  induced   tree  $T_1$  of  $G$. By  definition  of a  maximum  induced  tree,  $T_1$   contains all but a vertex $w$ of $C$.  Let  $T_2$  be  the  component  of  $G$   containing  $w$    and all  the  trees  joined  to  $G$ at  $w$. Then  $T_2$  is a tree  since  is a simple connected  graph  without  cycles.  Since all  the   vertices  of   $G$  are   in   $T_1\cup T_2$,    $\{T_1,T_2\}$  is a minimal tree cover of  $G$. As a consequence  $T(G)=2$.
\hfill $\Box$
\begin{proposition}\label{prop2cactus}
Let $G(V,E)$  be a cactus  graph with  at  least two  cycles.  Then  the   tree  cover  number $T(G)\ge 3$.
\end{proposition}

\noindent\MakeUppercase{Proof}:

Suppose  $G$ has  exactly  two cycles  $C_1$ and $C_2$. If  $v$ is a  cut  vertex of $C_1$ and $C_2$ then $G-v$ is a union  of two trees  $T_1$ and $T_2$. Hence vertices  of $G$  are  covered  by  induced trees  $T_1, T_2$ and $\{v\}$. If  there is a path  connecting $u\in C_1$  and $v\in C_2$,  then  removing  a  vertex  $w_1\in C_1$  such that  $w_1\ne u$ and  another  vertex  $w_2\in C_2$ such  that $w_2\ne v$  produces  three  induced trees $T_1,T_2,T_3$ such  that   all  vertices  of  $G$ are covered  by  these  three trees. Hence $T(G)\ge 3$ in the case  $G$  has exactly  two  cycles.

Suppose  $G$  has   three or more  cycles. Then  an  induced subgraph consisting  of  exactly  two  cycles  has  tree cover number of  at least   three. Hence $T(G)\ge 3$.
\hfill $\Box$
\begin{proposition}\label{prop3cactus}
Let $G(V,E)$  be a cactus  graph with  at  least two  cycles.  Then  $G$  satisfies  $\gccp$.
\end{proposition}

\noindent\MakeUppercase{Proof}:

Let $G(V,E)$  be a cactus  graph with  at  least two  cycles.  Since  $G$  is  outerplanar  we know from \cite{FS} that  $\msr(G) = |G| - T(G)$.  From proposition  \ref{prop2cactus} we  get  $T(G)\ge 3$  and  from  the Corollary \ref{coro1cactus}   we get  $\msr (\overline{G})\le 5$. Therefore $|G|-T(G)+5\le |G|+2$.
\hfill $\Box$

From  the results  above  we  get  the next result.
\begin{theorem}\label{teocactus}
Let $G(V,E)$  be a cactus  graph. Then  $G$  satisfies  $\gccp$.
\end{theorem}

\noindent\MakeUppercase{Proof}:

Let $G(V,E)$  be a cactus  graph.  Then  if   $G$ has no  cycles,   then  $G$  is a tree  we  know  from  \cite{FW2} that $\msr(\overline{G})\le 3$ and $G$ satisfies $\gccp$. If   $G$   is a unicyclic  graph then    from \cite{LH2} we  get  that  $G$ satisfies  $\gccp$ .  Finally,  if   $G$ has  at least  two  cycles   then  from  proposition \ref{prop3cactus} $G$  satisfies $\gccp$.
\hfill $\Box$

\begin{example} The  cactus  graph $G$ shown in figure \ref{cactus1} can be  labeled  in  such  a way   that $G$  satisfies  the definition  of  C-$\delta$ graph. From  the above  theorem $G$  satisfies $\gccp$.
\begin{center}
\includegraphics[height=60mm]{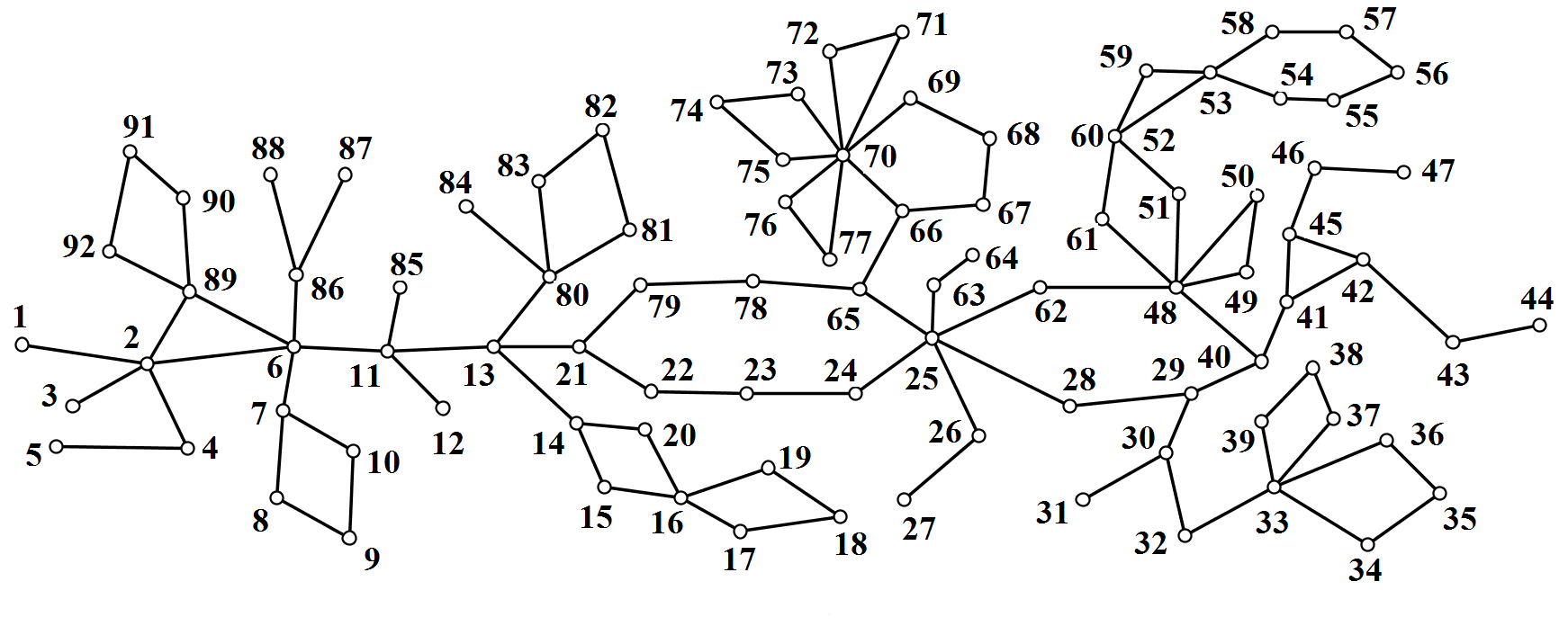}
\end{center}
\vspace{-0.1in}\begin{figure}[h]
\centering
\caption{Cactus Graph } \label{cactus2}
\end{figure}
\end{example}

\section{Conclusion}

Proving   GCC$_+$  for  cactus  graph  give us  other  way  to  prove  the result  for  several  families   like  trees, cycles,  chains of  cycles,  uncyclic graphs,  and   others    which   were  proved   by using  combinatorial approach ,  coloring and  other  techniques.  The  way  used  in  this paper  could  be used  to prove  GCC$_+$  for  several others   infinite  families  of simple connected  graphs  and  a  good approach  in proving GCC$_+$.
\section{Acknowledment}
I  would like  to  thanks   to  my  advisor  Dr. Sivaram Narayan  for   his  guidance and  suggestions  of  this  research.  Also  I  want  to  thank  to  the math  department  of  University  of  Costa  Rica and  Universidad  Nacional Estatal  a Distancia  because  their  sponsorship  during my  dissertation  research  and    specially  thanks  to   the math  department    of  Central Michigan University   where  I  did  the  researh   for  this  paper.
\input{bibliography}

\end{document}

%% file: bibliography.tex
\renewcommand{\baselinestretch}{1} 